\documentclass[reqno, english, 12pt,final]{amsart}
\def\vint_#1{\mathchoice%
          {\mathop{\kern 0.2em\vrule width 0.6em height 0.69678ex depth -0.58065ex
                  \kern -0.8em \intop}\nolimits_{\kern -0.4em#1}}%
          {\mathop{\kern 0.1em\vrule width 0.5em height 0.69678ex depth -0.60387ex
                  \kern -0.6em \intop}\nolimits_{#1}}%
          {\mathop{\kern 0.1em\vrule width 0.5em height 0.69678ex
              depth -0.60387ex
                  \kern -0.6em \intop}\nolimits_{#1}}%
          {\mathop{\kern 0.1em\vrule width 0.5em height 0.69678ex depth -0.60387ex
                  \kern -0.6em \intop}\nolimits_{#1}}}
\def\vintslides_#1{\mathchoice%
          {\mathop{\kern 0.1em\vrule width 0.5em height 0.697ex depth -0.581ex
                  \kern -0.6em \intop}\nolimits_{\kern -0.4em#1}}%
          {\mathop{\kern 0.1em\vrule width 0.3em height 0.697ex depth -0.604ex
                  \kern -0.4em \intop}\nolimits_{#1}}%
          {\mathop{\kern 0.1em\vrule width 0.3em height 0.697ex depth -0.604ex
                  \kern -0.4em \intop}\nolimits_{#1}}%
          {\mathop{\kern 0.1em\vrule width 0.3em height 0.697ex depth -0.604ex
                  \kern -0.4em \intop}\nolimits_{#1}}}

\usepackage{a4wide}
\usepackage{mathtools}
\usepackage{amsthm}
\usepackage{amssymb}
\usepackage{hyperref}
\usepackage{color}
\usepackage[obeyFinal]{todonotes}
\mathtoolsset{showonlyrefs}
\usepackage{enumerate}
\usepackage[obeyFinal]{todonotes}

\usepackage{csquotes}

\newcommand{\R}{\mathbb{R}}
\newcommand{\N}{\mathbb{N}}

\renewcommand{\P}{\mathcal{P}}

\renewcommand{\L}{\mathcal{L}}
\newcommand{\A}{\mathcal{A}}
\newcommand{\og}{\mathrm{OptGeo}}
\newcommand{\opt}[2]{\mathrm{Opt}(#1,#2)}
\newcommand{\id}{\mathrm{id}}
\newcommand{\m}{\mathfrak{m}}
\newcommand{\geo}{\mathrm{Geo}}
\newcommand{\spt}{\mathrm{spt}\,}
\newcommand{\restr}{\mathrm{restr}}
\renewcommand{\d}{\,\mathrm{d}}
\newcommand{\ent}[1]{\mathrm{Ent_\infty}\left( #1 \right)}
\newcommand{\entn}[1]{\mathrm{Ent_N}\left(#1\right)}
\newcommand{\Ent}{\mathrm{Ent}}

\newcommand{\ymt}{\frac{(t_3-t_2)}{(t_3-t_1)}}
\newcommand{\yt}{\frac{(t_2-t_1)}{(t_3-t_1)}}
\newcommand{\DC}{\mathcal{DC}}
\renewcommand{\liminf}{\varliminf}

\newtheorem{theorem}{Theorem}[section]
\newtheorem{lemma}[theorem]{Lemma}
\newtheorem{proposition}[theorem]{Proposition}
\newtheorem{corollary}[theorem]{Corollary}
\theoremstyle{definition}
\newtheorem{definition}[theorem]{Definition}

\theoremstyle{remark}
\newtheorem{remark}[theorem]{Remark}
\begin{document}
\title{Equivalent definitions of very strict $CD(K,N)$ -spaces}
\author{Timo Schultz}
\address{University of Jyvaskyla
\\ Department of Mathematics and Statistics
\\ P.O.Box 35
\\ FI-40014 University of Jyvaskyla}
\email{timo.m.schultz@jyu.fi}
\date{\today}

\begin{abstract} We show the equivalence of the definitions of very strict $CD(K,N)$ -condition defined, on one hand, using (only) the entropy functionals, and on the other, the full displacement convexity class $\DC_N$. In particular, we show that assuming the convexity inequalities for the critical exponent implies it for all the greater exponents. We also establish the existence of optimal transport maps in very strict $CD(K,N)$ -spaces with finite $N$.
\end{abstract}
\maketitle
\todo[inline]{{\color{red} Ideas/remarks: (1) By assume ``pointwise strict CD'' -condition (i.e. convexity for $\rho^F_t$), one could probably prove existence of map in the following way. Find $\pi^1$ and $\pi_2$ so that $\pi^1_t=\pi^2_t$ and $\pi^1_i\perp \pi^2_i$ for $i=0,1$. Then, by looking at $\lambda\pi^1+(1-\lambda)\pi^2$, and letting $\lambda\to\infty$ one gets contradiction with the pointwise convexity.

(2) Using pointwise convexity, one might be able to get that $ local vsCD + global CD \Rightarrow global vsCD$.   }}
\todo[inline]{Tarkista onko convexity inequalityyn viittaaminen riittävän selkeetä (eli että ei tule sekaannusta siita, toimiiko se tietylle siirrolle vai kaikille painotetuille)}
\section{Introduction}
Synthetic notions of curvature (bounds) have established their position in geometric analysis both as a tool to study geometric and analytic properties of non-smooth spaces, and as a new approach to attack problems even in the smooth setting. The framework present in this paper is the generalisation of Ricci curvature lower bounds to metric measure spaces, more precisely the setting of $CD(K,N)$-spaces introduced in the seminal papers of Lott--Villani \cite{Lott--Villani} and Sturm \cite{Sturm, Sturm2} based on a concept of displacement convexity of certain entropy functionals introduced by McCann \cite{McCann1997}. \todo{pitaako viitata Cordero-Erasquin, McCann, Schmuckenschlägeriin ja von Renesse-Sturmiin?}

The definitions of Sturm and Lott--Villani of $CD(K,N)$-spaces both share two notable properties, namely they are true generalisations of the notion of Ricci curvature lower bounds of (weighted) Riemannian manifolds, and, keeping in mind the Gromov's precompactness theorem for Riemannian manifolds sharing a common Ricci lower bound, they are stable under suitable convergence of metric measure spaces. The definitions of $CD(K,N)$-spaces by Sturm, and by Lott and Villani are different, but under an additional (essential) non-branching assumption of the spaces in question, these two notions of $CD(K,N)$-spaces agree. However, the non-branching property, while giving many desired results for $CD(K,N)$-spaces \cite{Sturm2, Cavalletti-Milman, Gigli-Rajala-Sturm, Rajala2016, Cavalletti-Sturm, Cavalletti-Mondino2017, Ohta}, is not stable under any reasonable convergence even when coupled with the $CD(K,N)$-condition.

In this paper, we study convexity properties  of a pointwise density of transport plans in (possibly) branching $CD(K,N)$-spaces giving an equivalent definition (Proposition \ref{LOCALISATION}) for the so-called \emph{very strict $CD(K,N)$} -condition introduced in \cite{Schultz2018} (see also \cite{Ambrosio-Gigli-Savare}), analogous to the known characterisation of essentially non-branching $CD(K,N)$-spaces, see \cite{Cavalletti-Milman}. Having the pointwise definition in hand, we prove Theorem \ref{equiviCD}, the equivalence of very strict $CD(K,N)$ -condition and its Lott--Villani type analogue (see Section \ref{PreliCD} for the precise definitions). 

The main difference in the definitions by Sturm and by Lott--Villani is that while Sturm requires convexity to hold only for certain specific entropy functionals, namely the R\'enyi entropies, Lott and Villani require it to hold for \emph{all} functionals in the so-called displacement convexity class. Using the defining convexity properties of the functionals in the displacement convexity class, we deduce easily the equivalence of the two definitions of very strict $CD(K,N)$ -spaces from the pointwise convexity inequality.

To obtain the pointwise condition, we use Theorem \ref{exmap}, the existence of optimal transport maps between two measures absolutely continuous with respect to the reference measure proven in \cite{Schultz2018} in the infinite dimensional case. For completeness, we present here the proof in the finite dimensional case. In fact, we need a bit more than just the existence of transport map. We need the plan to be given by a map not only from the endpoints, but also from the intermediate points. 

As a byproduct, we prove Theorem \ref{exmaps2}, the existence of optimal transport map from a (boundedly supported) absolutely continuous measure to a singular one. We construct the plan given by a map by gluing together plans obtained between (absolutely continuous) intermediate points of the endpoints. We prove, in similar fashion to what is done in \cite{Rajala2013}, that the resulting plan satisfies the convexity inequality of reduced curvature dimension condition between any three points of the unit interval.

\subsection*{Acknowledgements}The author would like to thank Enrico Pasqualetto for suggestions and discussions that led to the present paper. The author also acknowledges the support by the Academy of Finland, projects \#314789 and \#312488.
\section{Preliminaries} \label{Prelit}
Standing assumptions of this paper for a metric measure space $(X,d,\m)$ are completeness and separability for the metric $d$, and local finiteness for the Borel measure $\m$. 

A metric space $(X,d)$ is said to be a length space, if the distance between any two points $x$ and $y$ is obtained by infimising the length of curves connecting $x$ and $y$. A constant speed curve parametrised on the unit interval with length equal to the distance between the endpoints is called a (constant speed) geodesic. The set of all constant speed geodesics endowed with the supremum metric is denoted by $\geo(X)$.

\subsection{Optimal mass transportation} We consider the Monge--Kantorovich formulation of the optimal transport problem with quadratic cost. Denote by $\P(X)$ the set of all Borel probability measures on $X$. We define the Wasserstein 2-distance $W_2$ between two Borel probability measures $\mu,\nu\in\P(X)$ as the infimum
\begin{align} W_2(\mu,\nu)\coloneqq \left( \inf_{\sigma\in{\A(\mu,\nu)}}\int_{X\times X}d^2(x,y)\d\sigma(x,y)\right)^{\frac12},
\end{align}
where $\A(\mu,\nu)\coloneqq \{\sigma\in\P(X\times X): P^1_\# \sigma=\mu, P_\#^2 \sigma=\nu\}$ is the set of admissible transport plans between $\mu$ and $\nu$. The existence of an admissible plan that realises the infimum is true in rather general setting, including ours \cite{Villani}.  Such a minimising admissible plan is called an optimal plan, and the set of optimal plans between measures $\mu$ and $\nu$ is denoted by $\opt{\mu}{\nu}$.

Denote by $\P_2(X)$ the set of all Borel probability measures with finite second moment, that is, those $\mu\in\P(X)$ which are of finite $W_2$-distance from a Dirac mass. Moreover, denote by $\P_2^{ac}(X)$ a further subset of $\P_2(X)$ of measures absolutely continuous with respect to the reference measure $\m$. 

We recall, that the Wasserstein distance $W_2$ defines an actual metric on the set $\P_2(X)$. The space $(\P_2(X),W_2)$ inherits also some properties from the base space $X$, namely the space $(\P_2(X),W_2)$ is complete and separable length space, if $(X,d)$ is. In the case of length spaces, we have the following useful characterisation of Wasserstein geodesics. A curve $t\mapsto \mu_t\in\P_2(X)$ is geodesic, if and only if there exists a measure $\pi\in\P(\geo(X))$ so that $(e_0,e_1)_\#\pi\in\opt{\mu_0}{\mu_1}$, and $\mu_t=(e_t)_\#\pi$ for all $t\in[0,1]$, where $\gamma\mapsto e_t(\gamma)\coloneqq \gamma_t$ is the evaluation map \cite{Lisini}. Such a probability measure $\pi$ is called optimal dynamical plan, or just optimal plan for short, and the set of all optimal dynamical plans from $\mu_0$ to $\mu_1$ is denoted by $\og(\mu_0,\mu_1)$. 

Recall, that for $\pi\in\og(\mu_0,\mu_1)$, we have that $(\restr_{t_1}^{t_2})_\#(F\pi)$ is still an optimal plan for all $t_1,t_2\in[0,1]$, $t_1<t_2$, and for all $F$ with $\int F\d\pi=1$, where $\restr_{t_1}^{t_2}\colon \geo(X)\to\geo(X)$, $(\restr_{t_1}^{t_2})(\gamma)(t)=\gamma(tt_2+(1-t)t_1)$.  For $\pi\in\P(\geo(X))$, we denote by $\pi^{-1}$ the pushforward measure of $\pi$ under the map $\gamma\mapsto \gamma^{-1}$, $\gamma^{-1}(t)\coloneqq\gamma(1-t)$.

\subsection{Synthetic Ricci curvature lower bounds} \label{PreliCD}
Based on the notion of displacement convexity, introduced by McCann \cite{McCann1997}, of suitable entropy functionals,  Sturm \cite{Sturm}, and independetly Lott and Villani \cite{Lott--Villani} introduced notions of Ricci curvature lower bounds for general (non-smooth) metric measure spaces.

We recall the definition of a more restrictive version of curvature dimension condition -- the so-called \emph{very strict CD(K,N)} -condition -- and, motivated by the existence result for optimal maps in the context of such spaces, we introduce a Lott--Villani type analogue of the very strict $CD(K,N)$ -condition.

For the definitions, we need to introduce some auxiliary notation. As building blocks, we define,  for $K\in\R$ and $N\in (0,\infty]$, coefficients $[0,1]\times \R_+\to\R\cup\{\infty\}$, $(t,\theta)\mapsto\sigma_{K,N}^{(t)}(\theta)$ as
\begin{align}
\label{reduced coefficients} \sigma_{K,N}^{(t)}(\theta)\coloneqq  \left\{\begin{array}{ll}t, & \mathrm{if\ } N=\infty 
\\\infty, &\mathrm{if\ }K\theta^2\ge N\pi^2\\
\frac{\sin(t\theta\sqrt{\frac{K}{N}})}{\sin(\theta\sqrt{\frac{K}{N}})},& \mathrm{if\ } 0<K\theta^2<N\pi^2
\\ t, &\mathrm{if\ } K=0
\\ \frac{\sinh(t\theta\sqrt{\frac{-K}{N}})}{\sinh(\theta\sqrt{\frac{-K}{N}})}, & \mathrm{if\ } K<0.
\\\end{array}\right.
\end{align}
Using these coefficients we further define, for $N\in(1,\infty]$, coefficients $\beta_{K,N}^{(t)}(\theta)$ and $\tau_{K,N}^{(t)}(\theta)$ as
\begin{align}&\beta_{K,N}^{(t)}(\theta)\coloneqq t^{1-N}(\sigma_{K,N-1}^{(t)}(\theta))^{N-1},\quad \mathrm{and}
\\ &\tau_{K,N}^{(t)}(\theta)\coloneqq t^{\frac1N}(\sigma_{K,N-1}^{(t)}(\theta))^{\frac{N-1}{N}}.
\end{align}

To be precise, we define for $t>0, N> 1$
\begin{align}\beta_{K,N}^{(t)}(\theta)\coloneqq \left\{\begin{array}{ll} e^{\frac{K}6(1-t^2)\theta^2}, &\mathrm{if\ } N=\infty
\\ \infty, &\mathrm{if\ }N<\infty,\ K\theta^2> (N-1)\pi^2\\
\left(\frac{\sin(t\theta\sqrt{\frac{K}{N-1}})}{t\sin(\theta\sqrt{\frac{K}{N-1}})}\right)^{N-1},& \mathrm{if\ } 0< K\theta^2\le (N-1)\pi^2
\\ 1,& \mathrm{if\ } N<\infty,\ K=0 
\\ \left(\frac{\sinh(t\theta\sqrt{\frac{-K}{N-1}})}{t\sinh(\theta\sqrt{\frac{-K}{N-1}})}\right)^{N-1}, & \mathrm{if\ } N<\infty,\ K<0,\end{array}\right.
\end{align}
and $\beta_{K,N}^{(0)}\equiv 1$.

For $N\in(1,\infty]$, define the entropy functionals $\Ent_N\colon \P_2(X)\to \R\cup\{\pm\infty\}$ as 
\begin{align}\Ent_N(\mu)\coloneqq -\int \rho^{-\frac1N}\d\mu,
\end{align}
for $N<\infty$, and 
\begin{align}\Ent_\infty(\mu)\coloneqq \int\log\rho\d\mu+\int\infty\d\mu^\perp.
\end{align}
Here $\mu=\rho\m+\mu^\perp$ with $\mu^\perp \perp\m$, and $\mu^\perp(\{\rho>0\})=0$. Further, for transport plan $\pi\in\P(\geo(X))$ with $(e_0)_\#\pi=\mu_0\in\P_2(X)$, and for $t\in[0,1],\ K\in\R$, define the distorted entropy
\begin{align}\Ent^{(t)}_{N,\pi}(\mu_0)\coloneqq -\int \left(\beta_{K,N}^{(t)}(d(\gamma_0,\gamma_1))\right)^\frac1N\rho_0(\gamma_0)^{-\frac1N}\d\pi(\gamma),
\end{align}
for $N<\infty$, and 
\begin{align} \Ent^{(t)}_{\infty,\pi}(\mu_0)\coloneqq \int \log\left(\frac{\rho_0(\gamma_0)}{\beta^{(t)}_{K,\infty}(d(\gamma_0,\gamma_1))}\right)\d\pi(\gamma)+\int\infty\d\mu_0^\perp.
\end{align}
\begin{definition}\label{Sturm vsCD}We say that a metric measure space $(X,d,\m)$ is a very strict $CD(K,N)$ -space, if for all $\mu_0,\mu_1\in\P_2^{ac}(X)$ with bounded supports, there exists $\pi\in\og(\mu_0,\mu_1)$ such that for all non-negative and bounded Borel functions $F\colon \geo(X)\to \R$ with $\int F\d\pi=1$, and for all $t_1, t_2\in[0,1]$, $t_1<t_2$, we have
\begin{align}\label{K-convexity} \entn{\tilde\mu_{t}}\le (1-t)\Ent^{(1-t)}_{N,\tilde\pi}(\tilde\mu_{0})+t\Ent^{(t)}_{N,\tilde\pi^{-1}}(\tilde\mu_{1})
\end{align} 
for all $t\in[0,1]$, where $\tilde\mu_t\coloneqq (e_t)_\#\tilde\pi\coloneqq (e_t)_\#(\restr_{t_1}^{t_2})_\#F\pi$.
\end{definition}

\begin{remark}The definition would make sense also without the assumption on the boundedness of the supports. In that case, the functionals $\Ent_\infty$ and $\Ent_{\infty,\pi}$ are not a priori well-defined for all $\mu\in\P_2(X)$, due to the fact that $\int (\rho\log\rho)_-\d\m$ might be $-\infty$. However, after requiring \eqref{K-convexity} to hold for $\mu_i$, $i\in\{0,1\}$, with $(\rho_i\log\rho_i)_+\in\mathrm{L}^1(\m)$, we know by \cite[Theorem 4.24]{Sturm}, that (for fixed $x_0\in X$) $\m(B(x_0,r))\le Ae^{(Br^2)}$ holds for all $r>1$, and thus $(\rho\log \rho)_-\in \mathrm{L}^1(\m)$ for all $\mu=\rho\m\in\P_2(X)$, see \cite{Ambrosio-Gigli-Savare}. 
\end{remark}

We will also use the definition of \emph{very strict $CD^*(K,N)$ -condition}, which one gets by modifying the above definitions (see \cite{Bacher-Sturm} for the definition of reduced curvature dimension condition). More precisely, one replaces the convexity inequality $\eqref{K-convexity}$ by  inequality
\begin{align}\label{K-*convexity} \Ent_N(\tilde\mu_t)\le -\int \sigma_{K,N}^{(1-t)}(d(\gamma_0,\gamma_1))\rho_0^{-\frac1N}(\gamma_0)+\sigma_{K,N}^{(t)}(d(\gamma_0,\gamma_1))\rho_1^{-\frac1N}(\gamma_1)\d\tilde\pi.
\end{align}

Definition \ref{Sturm vsCD} is a (possibly) more restrictive version of the strict $CD(K,N)$ -condition introduced in \cite{Ambrosio-Gigli-Savare}, and is given in the spirit of Sturm's original definition for curvature dimension condition. To define Lott--Villani type analogue of the condition, we need to introduce the so-called \emph{displacement convexity classes}, introduced by McCann in \cite{McCann1997}. 

We say, that a continuous and convex function $U\colon \R_+\to \R$ is in the displacement convexity class $\DC_N$ (of dimension $N\in(1,\infty]$), if $U(0)=0$, and if the function $s\mapsto u(s)$ is convex, where $u$ is defined as
\begin{align}u\colon(0,\infty)\to \R,\quad s\mapsto s^NU(s^{-N}), 
\end{align}
if $N<\infty$, and
\begin{align} u\colon\R\to\R,\quad s\mapsto e^sU(e^{-s}),
\end{align}
if $N=\infty$.

\begin{remark}
We recall, that the displacement classes are nested. Indeed, if $N<N'$, we have that $\DC_{N'}\subset\DC_N$. This can be seen for example by writing 
\[u_N(s)\coloneqq s^NU(s^{-N})=(s^{\frac{N}{N'}})^{N'}U((s^{\frac{N}{N'}})^{-N'})\eqqcolon u_{N'}(s^{\frac{N}{N'}})\]
as a composition of a convex and decreasing function $u_{N'}$ and concave function $s\mapsto s^{\frac{N}{N'}}$. If $N'=\infty$, one writes
\[u_N(s)=e^{N\log s}U(e^{-N\log s}),\]
and concludes again, by concavity of $s\mapsto \log s$, that $u_N$ is convex.
\end{remark}

For $U\in\DC_N$, define the (entropy) functional $U_\m\colon \P_2(X)\to \R\cup\{\infty\}$ as
\begin{align}U_\m(\mu)\coloneqq \int U\circ\rho\d\m+\int U'(\infty)\d\mu^\perp,
\end{align}
where $U'(\infty)\coloneqq \underset{s\to\infty}{\lim}\frac{U(s)}{s}\in\R\cup\{\infty\}$. Furthermore, for $\pi\in\P(\geo(X))$, $K\in\R$ and $t\in[0,1]$, define the functional $U_{\pi,\m}^{(t)}\colon \P_2(X)\to \R\cup\{\infty\}$ as
\begin{align}U_{\pi,\m}^{(t)}(\mu)&\coloneqq \int_{X}\int_{\geo(X)} U\left(\frac{\rho(\gamma_0)}{\beta_{K,N}^{(t)}(\gamma_0,\gamma_1)}\right) \beta_{K,N}^{(t)}(\gamma_0,\gamma_1)\d\pi_x(\gamma)\d\m(x)
\\&\phantom{=}+\int_X U'(\infty)\d\mu^\perp,
\end{align}
where $\{\pi_x\}$ is a disintegration of $\pi$ with respect to the evaluation map $e_0$. 
\begin{remark}
The functional $U_{\pi,\m}^{(t)}$ is not well-defined in general due to the non-uniqueness of the disintegration. However, the definition will be used only for $\pi\in\og(\mu,\nu)$, in which case the disintegration is unique up to $\mu$-measure zero set. Another cause of being ill-defined is the possible integrability issue, which may appear both for the positive and for the negative part  of $U\circ\rho$ (and $\beta U( \rho/\beta)$), creating $\infty-\infty$ situations. This can be seen by taking $U(s)=s\log s-s^{1-\frac1N}$ in the hyperbolic space. Because of these issues, we will use the above definitions only for measures with bounded support, in which case the functionals are well-defined, see e.g. \cite[Theorem 17.28]{Villani} for the proof.
\end{remark}

\begin{definition}\label{Lott--Villani vsCD} A metric measure space is said to satisfy the very strict $CD(K,N)$ condition \emph{in the spirit of Lott--Villani}, if for all $\mu_0,\mu_1\in\P_2^{ac}(X)$ with bounded supports, there exists $\pi\in\og(\mu_0,\mu_1)$ such that for all bounded non-negative Borel functions $F\colon \geo(X)\to \R$ with $\int F\d\pi=1$, and for all $t_1, t_2\in[0,1]$, $t_1<t_2$, we have
\begin{align} U_\m({\tilde\mu_{t}})\le (1-t)U^{(1-t)}_{\tilde\pi,\m}(\tilde\mu_{0})+tU^{(t)}_{\tilde\pi^{-1},\m}(\tilde\mu_{1})
\end{align} 
for all $t\in[0,1]$ and for all $U\in\DC_N$, where $\tilde\mu_t\coloneqq (e_t)_\#\tilde\pi\coloneqq (e_t)_\#(\restr_{t_1}^{t_2})_\#F\pi$.
\end{definition}

\begin{remark}By choosing $U_N(s)= -s^{1-\frac1N}$, for $N<\infty$, and $U_\infty(s)=s\log s$, one immediately sees that spaces satisfying Definition \ref{Lott--Villani vsCD} also satisfy Definition \ref{Sturm vsCD}.
\end{remark}

 \section{Existence of optimal maps}\label{existence of maps}
 In proving our main results in Section \ref{ekvi}, we will use the fact that the plan given by the definition of very strict $CD(K,N)$ -spaces is induced by a map. The case $N=\infty$ is covered in \cite{Schultz2018}, and the proof of the finite dimensional case follows along the same lines. For completeness, we will outline the proof of the finite dimensional case here. It should be pointed out, that with our definition of very strict $CD(K,N)$ -spaces, we do not a priori know that very strict $CD(K,N)$ -condition for finite $N$ implies the very strict $CD(K,\infty)$ -condition.
\begin{theorem}[Existence of optimal maps]\label{exmap}
Let $(X,d,\m)$ be a very strict $CD^*(K,N)$ $(CD(K,N))$ -space, and let $\mu_0,\mu_1\in \P_2^{ac}(X)$ with bounded supports. Let $\pi\in\og(\mu_0,\mu_1)$ be the optimal plan given by the very strict $CD^*(K,N)$ $(CD(K,N))$ -condition. Then $\pi$ is induced by a Borel map $T\colon X\to \geo(X)$, i.e. $\pi=T_\#\mu_0$ with $e_0\circ T=\id$.
\end{theorem}

\begin{remark} If we remove in Definition \ref{Sturm vsCD} the assumption of the boundedness of the supports of $\mu_0$ and $\mu_1$, we may remove it also from Theorem \ref{exmap}.
\end{remark}

\begin{proof}
Let $N<\infty$, and $\mu_0,\mu_1\in\P_2^{ac}(X)$.  Furthermore, let $\pi\in\og(\mu_0,\mu_1)$ be the optimal plan given by  the definition of very strict $CD^*(K,N)$ -space. Suppose that $\pi$ is not induced by a map. Towards a contradiction, we will show that there exist plans $\pi^1,\pi^2\ll \pi$, and times $t_1$ and $t_2$ sufficiently close to each other so that $\mu^1_{t_i}=\mu^2_{t_i}$ and $\mu^1_{t_{i+1}}\perp\mu^2_{t_{i+1}}$. 

We begin by doing some reductions. First of all, by writing the whole space $X$ as a union of bounded sets, we may assume that the length of the geodesics in the support of $\pi$ is bounded by some constant $C$, and since $\spt\m$ is proper, we may also assume that $\spt\mu_0$ is compact. Furthermore, by dividing the interval $[0,1]$ into sufficiently small subintervals $I_{j}$, and looking at the restriction measures $(\restr_{I_j})_\#\pi$, we may assume that 
	\begin{align}\label{eq:localisation}\sigma_{K,N}^{(t)}(\theta)\in [(1-\varepsilon)t, (1-\varepsilon)^{-1}t]
	\end{align} 
for all $t\in[0,1]$ and $\theta\le C$. Here $\varepsilon>0$ is chosen so that $(1-\varepsilon)^42^{\frac1N}>1$.

Next, as was done in \cite{Schultz2018}, we find times $T,S\in(0,1)$, $T<S$, and optimal plans $\pi^1,\pi^2\ll \pi$ so that $\mu_T^1=\mu_T^2$ and $\mu_S^1\perp\mu_S^2$, where $\mu_t\coloneqq (e_t)_\#\pi$ for all $t\in[0,1]$. We refer to \cite{Schultz2018} for the arguments and the construction. Let then $n\in\N$ be such that 
	\begin{align}\label{TS}
	\frac{t}{t+\frac1n}\left(\frac{1-(t+\frac1n)}{1-t}\right)\ge \frac1{(1-\varepsilon)^4}2^{-\frac1N},
	\end{align}	
for $t\in[T,S]$. Again, by the arguments used in \cite{Schultz2018},  we find times $t_1,t_2\in[T,S]$, $t_1<t_2$, with $\lvert t_2-t_1\rvert<\frac1n$, and optimal plans $\bar\pi^1,\bar\pi^2$ such that $\bar\mu_{t_1}^1=\bar\mu_{t_2}^1$ and $\mu^1_{t_2}\perp \mu^2_{t_2}$. 

Now we are ready to arrive to a contradiction by similar computations as was done in \cite{Rajala2012}. We first use the convexity of the entropy along $\frac12(\bar\pi^1+\bar\pi^2)$ between points $0$, $t_1$ and $t_2$, then along $\bar\pi^1$ and $\bar\pi^2$ separately between points $t_1$, $t_2$ and $1$. Also the inequality \eqref{TS} is used both times with the convexity inequality. Then we use the bound \eqref{eq:localisation} and finally arrive to a contradiction.
	\begin{align}\int(\bar\rho^1_{t_1})^{1-\frac1N}\d\m&\ge(1-\varepsilon)^2 \frac{t_2-t_1}{t_2}2^{\frac1N-1}\int((\bar\rho^1_0)^{1-\frac1N}+(\bar\rho^2_0)^{1-\frac1N})\d\m
	\\ &\phantom{\ge}+(1-\varepsilon)^2\frac{t_1}{t_2}2^{\frac1N-1}\left(\int(\bar\rho^1_{t_2})^{1-\frac1N}\d\m+\int(\bar\rho^2_{t_2})^{1-\frac1N}\d\m\right)
	\\ &> (1-\varepsilon)^4\frac{t_1}{t_2}\frac{(1-t_2)}{(1-t_1)}2^{\frac1N}\int(\bar\rho_{t_1}^1)^{1-\frac1N}\d\m\ge\int(\bar\rho^1_{t_1})^{1-\frac1N}\d\m.
	\end{align}
Here $\bar\rho^i_t$ is the density of $(e_t)_\#\bar\pi^i$ with respect to $\m$. In the case of very strict $CD(K,N)$ -space, the proof is exactly the same after replacing $\sigma_{K,N}^{(t)}$ by $\tau_{K,N}^{(t)}$ in the condition \eqref{eq:localisation}.
\end{proof}

As a corollary, we get the existence of an optimal map from absolutely continuous measure to singular one, by  approaching the singular endpoint with absolutely continuous intermediate points. Combined with construction similar to the one used in \cite{Rajala2013}, we arrive to the following theorem.

\begin{theorem} \label{exmaps2}
Let $(X,d,\m)$ be a very strict $CD^*(K,N)$ -space with $N<\infty$, and $\mu_0\in\P_2^{ac}(X)$ and $\mu_1\in\P_2(X)$, $\spt\mu_1\subset \spt\m$, probability measures with bounded support. Then there exists $\pi\in\og(\mu_0,\mu_1)$ along which the convexity inequality \eqref{K-*convexity} holds between any points $t_1<t_2<t_3$ (with $\tilde\mu_t=\mu_t=(e_t)_\#\pi$) for the entropy $\Ent_N$. Moreover, $\pi$ is induced by a map from $\mu_0$.
\end{theorem}
\begin{remark}
We do not claim, that the convexity would hold along $F\pi$, where $F$ is arbitrary bounded non-negative Borel function with $\int F\d\pi=1$. In fact, the proof below will in some cases produce a geodesic $(\mu_t)$ such that for any lift $\pi$ of $(\mu_t)$ this is known to be false. 
\end{remark}
The idea of the proof of the above theorem is fairly simple. First of all, by approximating the possibly singular measure $\mu_1$ by absolutely continuous ones, one obtains a geodesic $\mu_t$ with $\mu_t\ll\m$ due to the lower semi-continuity of the entropy $\Ent_N$. Then, by compactness of midpoints, there exist $t$-intermediate points $\mu_t$, $t\in\{\frac12,\frac34,\frac78,\dots\}$, that are absolutely continuous, and minimise the entropy $\Ent_N$ among all midpoints of the previous point and $1$. Now taking $\pi^i\in\og(\mu_{\frac{2^i-1}{2^i}},\mu_{\frac{2^{i+1}-1}{2^{i+1}}})$ given by Theorem \ref{exmap}, and concatenating them, one obtains in the limit a plan with desired properties.

In the proof we will use the following lemma.
\begin{lemma}\label{abscontinter} Let $(X,d,\m)$ be a very strict $CD^*(K,N)$ -space with $N<\infty$, and $\mu_0\in\P_2^{ac}(X)$ and $\mu_1\in\P_2(X)$, $\spt\mu_1\subset \spt\m$, probability measures with bounded support, and let $\mu_\frac12$ be a midpoint of $\mu_0$ and $\mu_1$ minimising the entropy among all midpoints. Then $\mu_\frac12\in\P_2^{ac}(X)$.
\end{lemma}

\begin{proof} Clearly, we may assume that $K<0$. Let $\mu_1^i$ be a sequence of absolutely continuous measures converging to $\mu_1$, and having (uniformly) bounded support. Let $\pi_i\in\og(\mu_0,\mu_1^i)$ be a sequence satisfying the convexity inequality \eqref{K-*convexity} and (sub)converging to some $\pi\in\og(\mu_0,\mu_1)$. Then, by lower semi-continuity of the entropy, we have
\begin{align} \Ent_N(\mu_\frac12)&\le \liminf_i \Ent_N(\mu^i_\frac12)\le\liminf_i -\int\sigma_{K,N}^{\frac12}(d(\gamma_0,\gamma_1))\rho_0^{-\frac1N}(\gamma_0)\d\pi^i
\\ &\le -\sigma_{K,N}^\frac12(D)\Ent_N(\mu_0)<0,
\end{align}
where $D$ is a bound for the diameters of the supports. Thus, we know that $\mu_\frac12$ is not purely singular. Let now $\mu_\frac12=\rho_\frac12\m+\mu_\frac12^\perp$ be the Lebesgue decomposition of $\mu_\frac12$, and let $A$ be a Borel set on which $\mu_\frac12^\perp$ is concentrated and with $(\rho_\frac12\m)(A)=0$.  We want to show that $\mu_\frac12^\perp(A)=0$. Suppose that this is not the case, and define $\tilde\mu_j\coloneqq (e_j)_\#\pi\lvert_{e_\frac12^{-1}(A)}$ for $j\in\{0,1\}$. Since $\tilde \mu_0$ is absolutely continuous with respect to $\m$, there exists, by taking the minimiser of the entropy (which exists by compactness of midpoints), a midpoint $\tilde\mu_\frac12$ of $\tilde\mu_0$ and $\tilde\mu_1$ which is not purely singular. Hence, $\hat\mu_\frac12\coloneqq \rho_\frac12\m+\tilde\mu_\frac12$ is a midpoint of $\mu_0$ and $\mu_1$ with
\[\Ent_N(\hat\mu_\frac12)<\Ent_N(\mu_\frac12),\]
which contradicts the assumption of $\mu_\frac12$ realising the minimum of the entropy.
\end{proof}
\begin{proof}[Proof of Theorem \ref{exmaps2}]
Let $\mu_0\in\P_2^{ac}(X)$ and $\mu_1\in\P_2(X)$, $\spt\mu_1\subset \spt\m$, be probability measures with bounded support. 
Since the space is boundedly compact, we know that the set of midpoints $\mathcal{M}(\mu_0,\mu_1)$ of $\mu_0$ and $\mu_1$ is compact. Moreover, the entropy $\Ent_N$ is lower semi-continuous on  $\mathcal{M}(\mu_0,\mu_1)$ due to the finiteness of $\m$ on bounded sets. Thus, there exists a midpoint $\mu_{\frac12}\in\mathcal{M}(\mu_0,\mu_1)$ that minimises the entropy among midpoints. By induction we get a sequence of $t_i$-intermediate points $(\mu_{t_i})_{i\in\N}$, where $t_i=(1-2^{-i})$, and $\mu_{t_i}$ minimises the entropy among all midpoints of $\mu_{t_i}$ and $\mu_1$. 

By Lemma \ref{abscontinter} we have that $\mu_{t_i}\ll\m$. Thus, for each $i\in\N$, there exists $\pi_i\in\og(\mu_{t_{i-1}},\mu_{t_i})$ satisfying the very strict $CD^*(K,N)$ -condition, hence is induced by a map $T_i$ from $\mu_{t_{i-1}}$. Consider now the decreasing sequence $(A_i)$ of sets
\[A_i\coloneqq \{\pi\in\og(\mu_0,\mu_1): (\restr_{t_{j-1}}^{t_{j}})_\#\pi=\pi_j \mathrm{\ for\ all\ } j\le i\}.\]
We will show that the intersection $A\coloneqq \cap A_i$ is singleton, and that the unique element of $A$ satisfies the desired conditions. Since the sequence is nested, to show that $A$ is non-empty, it suffices to show that each $A_i$ is compact. Each $A_i$ is tight, since the set $(\restr_0^{\frac12})^{-1}(\spt\pi^1)$ is compact due to the continuity of the map $(\restr)_0^{\frac12}$ and Arzelà-Ascoli theorem. To see that $A_i$ is closed, take a converging sequence $\tilde\pi_n\in\A_i$, $\tilde\pi_n\longrightarrow\tilde\pi$. Then $(\restr_{t_j}^{t_{j+1}})_\#\tilde\pi_n\longrightarrow (\restr_{t_j}^{t_{j+1}})_\#\tilde\pi$, and hence $(\restr_{t_j}^{t_{j+1}})_\#\tilde\pi=\pi_j$. Therefore $A_i$ is compact, and $A$ is non-empty.

Let now $\pi\in A$. Then, for all $i\in\N$, $(\restr_{0}^{t_i})_\#\pi$ is induced by a map $T_i$ due to the fact that $(\restr_{t_{i-1}}^{t_i})_\#\pi$ is induced by a map (see, e.g.\ \cite[Lemma 4]{Schultz2018}). When $i<j$, we have that $T_i=\restr_0^{{t_i}/{t_j}}\circ T_j$. Thus, we have by completeness of $X$ that $T_{i}$ converges pointwise to some $T$. Indeed, for any $x\in X$, the sequence $T_i(x)$ is a Cauchy sequence. Hence, by dominated convergence we have for any continuous and bounded function $f\colon \geo({X}) \to \R$, that
\begin{align} \int f\d (T_i)_\#\mu_0= \int f\circ T_i\d\mu_0\longrightarrow \int f\circ T\d\mu_0=\int f\d T_\#\mu_0
\end{align} 
giving the weak convergence $(\restr_0^{t_i})_\#\pi\longrightarrow T_\#\mu_0$. On the other hand, we know that $(\restr_0^{t_i})_\#\pi\longrightarrow \pi $. Hence, the plan $\pi$ is induced by a map.

Let us now prove the convexity of the entropy along $\pi$. The steps are similar to the ones in \cite{Rajala2013}. We will first prove, that the convexity holds between points $\delta$, $\frac12$ and $1$, where $\delta$ is arbitrarily small. Let $\delta\in(0,\frac12)$. Suppose now that the claim is not true. Then there exists an interval $I=(a,b)\subset (0,\infty)$ so that 
\begin{align} \label{failure}
\int_{l^{-1}(I)}\rho_{\frac12}^{-\frac1N}(\gamma_{\frac12})\d\pi < \int_{l^{-1}(I)} \sigma_{K,N}^{\frac{\frac12}{1-\delta}}(d(\gamma_\delta,\gamma_1))\rho^{-\frac1N}_\delta(\gamma_\delta)+\sigma_{K,N}^{\frac{\frac12-\delta}{1-\delta}}(d(\gamma_\delta,\gamma_1))\rho^{-\frac1N}_1(\gamma_1)\d\pi,
\end{align} 
where $l\colon \geo(X)\to \R$ is the map sending a geodesic to its length. By continuity of the distortion coefficients we may assume, by subdividing the interval further, that 
\begin{align}\label{epsilondelta}(1-\varepsilon)\sigma^\alpha((1-\delta)a)\le \sigma^\alpha((1-\delta)b),\end{align} 
where $\alpha\in\{\frac{\frac12}{1-\delta},\frac{\frac12-\delta}{1-\delta}\}$, and $\varepsilon$ is chosen so that
\begin{align}\label{epsilonfailure}\begin{split}
\int_{l^{-1}(I)}\rho_{\frac12}^{-\frac1N}(\gamma_{\frac12})\d\pi&  
\\< (1-\varepsilon)&\int_{l^{-1}(I)} \sigma_{K,N}^{\frac{\frac12}{1-\delta}}(d(\gamma_\delta,\gamma_1))\rho^{-\frac1N}_\delta(\gamma_\delta)+\sigma_{K,N}^{\frac{\frac12-\delta}{1-\delta}}(d(\gamma_\delta,\gamma_1))\rho^{-\frac1N}_1(\gamma_1)\d\pi. \end{split}
\end{align} 

Let $\pi^I\coloneqq \pi\lvert_{l^{-1}(I)}$, and let $\mu_j^I=(e_j)_\#\pi^I$ for $j\in\{0,\delta,1\}$. Let $\mu^i_1\longrightarrow \mu_1^I$ be a sequence of absolutely continuous measures (with equibounded support) for which
\[\int (\rho_1^i)^{1-\frac1N}\d\m\longrightarrow  \int (\rho_1^I)^{1-\frac1N}\d\m.\]
This can be done simply by approximating separately the singular part of $\mu_1^I$, due to the lower semi-continuity of the entropy. Let now $\pi^i\in\og(\mu_0^I,\mu_1^i)$ be such that the converse of \eqref{failure} holds for $\pi^i$ between points $\delta_i$, $\frac12$ and $1$, where $\delta_i\longrightarrow \delta$ with $\tilde\mu_{\delta_i}=\mu_\delta^I$. Finally, define \[\tilde\pi^i\coloneqq \pi\lvert_{\geo(X)\setminus \l^{-1}(I)}+\pi^i.\]
We may assume, that $\tilde\pi^i\longrightarrow \tilde\pi\in\og(\mu_0,\mu_1)$ weakly. By c-cyclical monotonicity (see \cite[Proposition 1]{Rajala2013}), and by weak convergence,  we know that $\tilde\pi^i(l^{-1}(I))\longrightarrow 1$. Thus,

\begingroup
\allowdisplaybreaks
\begin{equation} \label{arvio1}\begin{split}\Ent_N(\tilde\mu_{\frac12})& \le\liminf_{i\to\infty}\Ent_N(\tilde\mu_1^i)= \liminf_{i\to\infty}\left[-\int_{\geo(X)\setminus l^{-1}(I)}(\tilde\rho_{\frac12}^i)^{-\frac1N}(\gamma_{\frac12})\d\pi-\int_{\geo(X)}(\tilde\rho_{\frac12}^i)^{-\frac1N}(\gamma_{\frac12})\d\pi^i\right]
\\ &\le\liminf_{i\to\infty}\left[ -\int_{\geo(X)\setminus l^{-1}(I)}\rho^{-\frac1N}_{\frac12}(\gamma_{\frac12})\d\pi-\int(\rho^i_{\frac12})^{-\frac1N}(\gamma_{\frac12})\d\pi^i\right]
\\ &\le \liminf_{i\to\infty}\left[-\int_{\geo(X)\setminus l^{-1}(I)}\rho^{-\frac1N}_{\frac12}(\gamma_{\frac12})\d\pi\right.
\\ &\phantom{=}\left.-\int\sigma^{\frac{\frac12}{1-\delta_i}}_{K,N}(d(\gamma_{\delta_i},\gamma_1))(\rho^i_{\delta_i})^{-\frac1N}(\gamma_{\delta_i})+\sigma^{\frac{\frac12-\delta_i}{1-\delta_i}}_{K,N}(d(\gamma_{\delta_i},\gamma_1))(\rho^i_{1})^{-\frac1N}(\gamma_{1})\d\pi^i\right]
\\ &\le \liminf_{i\to\infty}\Bigg[-\int_{\geo(X)\setminus l^{-1}(I)}\rho^{-\frac1N}_{\frac12}(\gamma_{\frac12})\d\pi
\\&\phantom{=}-\int_{\l^{-1}(I)}\sigma^{\frac{\frac12}{1-\delta_i}}_{K,N}(d(\gamma_{\delta_i},\gamma_1))(\rho^i_{\delta_i})^{-\frac1N}(\gamma_{\delta_i})\d\pi^i-\int_{\l^{-1}(I)}\sigma^{\frac{\frac12-\delta_i}{1-\delta_i}}_{K,N}(d(\gamma_{\delta_i},\gamma_1))(\rho^i_{1})^{-\frac1N}(\gamma_{1})\d\pi^i\Bigg] \end{split}
\end{equation}
\endgroup
due to the lower semi-continuity of the entropy, the fact $\tilde\rho_\frac12\le \rho_\frac12^i$ everywhere, $\tilde\rho_\frac12(\gamma_\frac12)\le\rho_\frac12(\gamma_\frac12)$ in $\geo(X)\setminus l^{-1}(I)$, and the convexity of the entropy along $\pi^i$. To arrive to a contradiction, we will need the following observation, which follows by the disintegration theorem, Hölder's inequality, and Jensen's inequality. \begin{align}\label{junktozero}\begin{split}
&\phantom{=}\int_{\geo(X)\setminus l^{-1}(I)}(\rho^i_{t})^{-\frac1N}(\gamma_{t})\d\pi^i=\iint \chi_{\geo(X)\setminus l^{-1}(I)}(\rho^i_{t})^{-\frac1N}\circ e_t\d\pi^i_x\d\mu_t^i(x)
\\&=\int(\rho^i_{t})^{-\frac1N}(x)\int \chi_{\geo(X)\setminus l^{-1}(I)}\d\pi^i_x\d\mu_t^i(x)=\int(\rho^i_{t})^{1-\frac1N}(x)\int \chi_{\geo(X)\setminus l^{-1}(I)}\d\pi^i_x\d\m(x)
\\&\le\left(\int \rho_t^i(x) \left(\int \chi_{\geo(X)\setminus l^{-1}(I)}\d\pi^i_x\right)^{\frac{N}{N-1}}\d\m(x)\right)^{\frac{N-1}{N}}(\m(\spt\mu_1^i))^{\frac1N}
\\&\le C\left(\int \rho_t^i(x) \int \chi_{\geo(X)\setminus l^{-1}(I)}\d\pi^i_x\d\m(x)\right)^{\frac{N-1}{N}}=C\pi_i(\geo(X)\setminus l^{-1}(I))\longrightarrow 0,\end{split}
\end{align}
when $i\longrightarrow \infty$. Hence, by \eqref{arvio1}
\begingroup
\allowdisplaybreaks
\begin{align}
&\Ent_N(\tilde\mu_{\frac12})
\\ &\le \liminf_{i\to\infty}\Bigg[-\int_{\geo(X)\setminus l^{-1}(I)}\rho^{-\frac1N}_{\frac12}(\gamma_{\frac12})\d\pi
\\&\phantom{=}-\sigma^{\frac{\frac12}{1-\delta_i}}_{K,N}((1-\delta_i)b)\int_{\l^{-1}(I)}(\rho^i_{\delta_i})^{-\frac1N}(\gamma_{\delta_i})\d\pi^i-\sigma^{\frac{\frac12-\delta_i}{1-\delta_i}}_{K,N}((1-\delta_i)b)\int_{l^{-1}(I)}(\rho^i_{1})^{-\frac1N}(\gamma_{1})\d\pi^i\Bigg]
\\ &= \liminf_{i\to\infty}\left[-\int_{\geo(X)\setminus l^{-1}(I)}\rho^{-\frac1N}_{\frac12}(\gamma_{\frac12})\d\pi+\sigma^{\frac{\frac12}{1-\delta_i}}_{K,N}((1-\delta_i)b)\Ent_N(\mu_\delta^I)+\sigma^{\frac{\frac12-\delta_i}{1-\delta_i}}_{K,N}((1-\delta_i)b)\Ent_N(\mu_1^i)\right.
\\&\phantom{=}\left.+\sigma^{\frac{\frac12}{1-\delta_i}}_{K,N}((1-\delta_i)b)\int_{\geo(X)\setminus l^{-1}(I)}(\rho^i_{\delta_i})^{-\frac1N}(\gamma_{\delta_i})\d\pi^i+\sigma^{\frac{\frac12-\delta_i}{1-\delta_i}}_{K,N}((1-\delta_i)b)\int_{\geo(X)\setminus l^{-1}(I)}(\rho^i_{1})^{-\frac1N}(\gamma_{1})\d\pi^i\right]
\\&\overset{\eqref{junktozero}}{=}-\int_{\geo(X)\setminus l^{-1}(I)}\rho^{-\frac1N}_{\frac12}(\gamma_{\frac12})\d\pi+\sigma^{\frac{\frac12}{1-\delta}}_{K,N}((1-\delta)b)\Ent_N(\mu_\delta^I)+\sigma^{\frac{\frac12-\delta}{1-\delta}}_{K,N}((1-\delta)b)\Ent_N(\mu_1^I)
\\&\overset{\eqref{epsilondelta}}{\le}-\int_{\geo(X)\setminus l^{-1}(I)}\rho^{-\frac1N}_{\frac12}(\gamma_{\frac12})\d\pi-(1-\varepsilon)\int_{l^{-1}(I)}\sigma^{\frac{\frac12}{1-\delta}}_{K,N}(d(\gamma_\delta,\gamma_1))\rho_\delta^{-\frac1N}(\gamma_\delta)+\sigma^{\frac{\frac12-\delta}{1-\delta}}_{K,N}(d(\gamma_\delta,\gamma_1))\rho_1^{-\frac1N}(\gamma_\delta)\d\pi
\\ &\overset{\eqref{epsilonfailure}}{<}-\int_{\geo(X)\setminus l^{-1}(I)}\rho^{-\frac1N}_{\frac12}(\gamma_{\frac12})\d\pi-\int_{l^{-1}(I)}\rho^{-\frac1N}_{\frac12}(\gamma_{\frac12})\d\pi=\Ent_N(\mu_\frac12),
\end{align}
\endgroup
which is a contradiction, since $\mu_\frac12$ was the minimiser of the entropy. Notice, that we actually proved a stronger version of the convexity between points $\delta$, $\frac12$ and $1$, namely that the convexity holds whenever the plan $\pi$ is restricted to $l^{-1}(I)$ for any open interval $I$.

To show that the convexity holds between $0$, $\frac12$ and $1$,  we use the convexity first between $\delta$, $\frac12$ and $1$, then between $0$, $\delta$, and $\frac12$, and then conclude by letting $\delta\longrightarrow 0$. Write 
\[[0,\infty)=\cup_{i\in\N}I_i,\]
where $I_i=[s_i,s_{i+1}]$ are intervals with equal length $\varepsilon>0$. Since the functions $\sigma^{(t)}_{K,N}$ are Lipschitz continuous with uniform Lipschitz constant $L$, we have 
\begingroup
\allowdisplaybreaks
\begin{align}
\Ent_N(\mu_{\frac12})&\le-\int \sigma^{\frac{\frac12}{1-\delta}}_{K,N}(d(\gamma_\delta,\gamma_1))\rho^{-\frac1N}_\delta(\gamma_\delta)+\sigma^{\frac{\frac12-\delta}{1-\delta}}_{K,N}(d(\gamma_\delta,\gamma_1))\rho^{-\frac1N}_1(\gamma_1)\d\pi
\\ &\le -\sum_{i}\sigma^{\frac{\frac12}{1-\delta}}_{K,N}((1-\delta)s_{i+1})\int_{l^{-1}(I_i)}\rho^{-\frac1N}_\delta(\gamma_\delta)\d\pi-\int \sigma^{\frac{\frac12-\delta}{1-\delta}}_{K,N}(d(\gamma_\delta,\gamma_1))\rho^{-\frac1N}_1(\gamma_1)\d\pi 
\\ &\le -\sum_{i}\sigma^{\frac{\frac12}{1-\delta}}_{K,N}((1-\delta)s_{i+1})\int_{l^{-1}(I_i)}\sigma^{\frac{\frac12-\delta}{\frac12}}_{K,N}(d(\gamma_0,\gamma_{\frac12}))\rho^{-\frac1N}_0(\gamma_0)\d\pi
\\ &\phantom{\le}-\sum_{i}\sigma^{\frac{\frac12}{1-\delta}}_{K,N}((1-\delta)s_{i+1})\int_{l^{-1}(I_i)}\sigma^{\frac{\delta}{\frac12}}_{K,N}(d(\gamma_0,\gamma_{\frac12}))\rho^{-\frac1N}_{\frac12}(\gamma_{\frac12})\d\pi
\\ &\phantom{=}-\int \sigma^{\frac{\frac12-\delta}{1-\delta}}_{K,N}(d(\gamma_\delta,\gamma_1))\rho^{-\frac1N}_1(\gamma_1)\d\pi
\\ &\le -\sum_{i}(1-L\lvert s_{i+1}-s_i\rvert)\sigma^{\frac{\frac12}{1-\delta}}_{K,N}((1-\delta)s_{i})\int_{l^{-1}(I_i)}\sigma^{\frac{\frac12-\delta}{\frac12}}_{K,N}(d(\gamma_0,\gamma_{\frac12}))\rho^{-\frac1N}_0(\gamma_0)\d\pi
\\ &\phantom{\le}-\sum_{i}(1-L\lvert s_{i+1}-s_i\rvert)\sigma^{\frac{\frac12}{1-\delta}}_{K,N}((1-\delta)s_{i})\int_{l^{-1}(I_i)}\sigma^{\frac{\delta}{\frac12}}_{K,N}(d(\gamma_0,\gamma_{\frac12}))\rho^{-\frac1N}_{\frac12}(\gamma_{\frac12})\d\pi
\\ &\phantom{=}-\int \sigma^{\frac{\frac12-\delta}{1-\delta}}_{K,N}(d(\gamma_\delta,\gamma_1))\rho^{-\frac1N}_1(\gamma_1)\d\pi
\\&\le -\sum_{i}(1-L\lvert s_{i+1}-s_i\rvert)\int_{l^{-1}(I_i)}\sigma^{\frac{\frac12}{1-\delta}}_{K,N}((1-\delta)d(\gamma_0,\gamma_1))\sigma^{\frac{\frac12-\delta}{\frac12}}_{K,N}(d(\gamma_0,\gamma_{\frac12}))\rho^{-\frac1N}_0(\gamma_0)\d\pi
\\ &\phantom{\le}-\sum_{i}(1-L\lvert s_{i+1}-s_i\rvert)\int_{l^{-1}(I_i)}\sigma^{\frac{\frac12}{1-\delta}}_{K,N}((1-\delta)d(\gamma_0,\gamma_1))\sigma^{\frac{\delta}{\frac12}}_{K,N}(d(\gamma_0,\gamma_{\frac12}))\rho^{-\frac1N}_{\frac12}(\gamma_{\frac12})\d\pi
\\ &\phantom{=}-\int \sigma^{\frac{\frac12-\delta}{1-\delta}}_{K,N}(d(\gamma_\delta,\gamma_1))\rho^{-\frac1N}_1(\gamma_1)\d\pi
\\ &\le (1-\varepsilon)\int\sigma^{\frac{\frac12}{1-\delta}}_{K,N}((1-\delta)d(\gamma_0,\gamma_1))\sigma^{\frac{\frac12-\delta}{\frac12}}_{K,N}(d(\gamma_0,\gamma_{\frac12}))\rho^{-\frac1N}_0(\gamma_0)+\sigma^{\frac{\frac12-\delta}{1-\delta}}_{K,N}(d(\gamma_\delta,\gamma_1))\rho^{-\frac1N}_1(\gamma_1)\d\pi
\\ &\phantom{=}(1-\varepsilon)\int \sigma^{\frac{\frac12}{1-\delta}}_{K,N}((1-\delta)d(\gamma_0,\gamma_1))\sigma^{\frac{\delta}{\frac12}}_{K,N}(d(\gamma_0,\gamma_{\frac12}))\rho^{-\frac1N}_{\frac12}(\gamma_{\frac12})\d\pi
\\ &\longrightarrow \int\sigma^{\frac12}_{K,N}(d(\gamma_0,\gamma_1))\rho^{-\frac1N}_0(\gamma_0)+\sigma^{\frac12}_{K,N}(d(\gamma_0,\gamma_1))\rho^{-\frac1N}_1(\gamma_1)\d\pi,
\end{align}
\endgroup
where we first let $\delta\longrightarrow 0$, and then $\varepsilon\longrightarrow 0$. In the first limit, we used dominated convergence with $C\rho^{1-\frac1N}_i$ as a dominant, and the explicit form of the distortion coefficients.

To show that the convexity holds between points $0$, $t$, and $1$, where $t\in(0,\frac12)$, one uses analogous computations as above, now using the convexity first between points $0$, $t$, and $\frac12$, and then between $\delta$, $\frac12$, and $1$, again letting $\delta\to 0$. 

Finally, the case for general $t\in(0,1)$ follows inductively -- after the observation that the convexity between $t_i$, $t_{i+1}$, and $1$ (and thus, between points $\frac12$, $t_i$, and $1$ by yet another induction) is of the stronger form, more precisely, the convexity holds when restricted to curves with length in an interval $[a,b]$ (converse inequality of \eqref{failure} with $\delta=t_i$, and $t_{i+1}$ in place of $\frac12$).

We have now shown that the convexity holds between points $0$, $t$, and $1$ for any $t\in(0,1)$. Next, we will turn into the proof of convexity between any three points $r<s<t$. It will follow analogously to the previous case after a couple of simple observations. First of all, if $r\in[t_i,t_{i+1})$, and $t\in[t_k,t_{k+1})$, $k>i$, then $\mu_{t_{i+1}}$ minimises the entropy among all $(t_{i+1}-r)/(t-r)$-intermediate points of $\mu_r$ and $\mu_t$. Furthermore, $\mu_j$ minimises the entropy among all $(t_j-t_{j-1})/(t-t_{j-1})$-intermediate points of $\mu_{t_{j-1}}$ and $\mu_t$ for all $j\in\{i+2,\dots, k\}$. The second observation needed is that the pushforward of a plan given by the definition of very strict $CD^*(K,N)$-space under the restriction map still satisfies the requirements of the very same definition. The only difference in the argument is that now instead of infinitely many steps in the induction argument, one only has a finite number of steps, and one special case, namely when $s\in(t_k,t)$. This special case, however, follows easily with the same arguments.
\end{proof}
 \section{Equivalent definitions of very strict $CD(K,N)$ -condition}\label{ekvi}
In this section we will prove that the definition of very strict $CD(K,N)$ -spaces is equivalent to an analogous pointwise convexity requirement for the density of a Wasserstein geodesic along optimal plan. This pointwise definition is then used to prove the equivalence of the definition of very strict $CD(K,N)$ and Lott--Villani-type analogous of the definition.

We will need the following simple lemma.
\begin{lemma}\label{imabs} Let $(X,d,\m)$ be a very strict $CD(K,N)$ -space, $\mu_0,\mu_1\in\P_2^{ac}(X)$ absolutely continuous measures with respect to the reference measure and with bounded supports, and $\pi\in\og(\mu_0,\mu_1)$ given by the definition of very strict $CD(K,N)$ -condition. Then $\mu_t\in\P_2^{ac}(X)$ for all $t\in(0,1)$.
\end{lemma}

\begin{proof}Suppose the claim is not true. Then there exists $\pi\in\og({\mu_0},{\mu_1})$ as in Definition \ref{Sturm vsCD} with $\mu_t\coloneqq (e_t)_\#\pi=\rho\m+\mu^\perp$, $\mu^\perp\perp\m$. Thus, there exists a Borel set $A\subset X$ so that $\mu^\perp(A)>0$ and $\m(A)=0$.  Let $\mathcal{A}\coloneqq e_t^{-1}(A)$, and define $\tilde\pi\coloneqq \pi\lvert_\mathcal{A}$. 

In the case $N=\infty$ we get a contradiction after restricting the plan $\pi$ further so that $\rho_0$ and $\rho_1$ are bounded, and hence the entropies $\Ent_\infty(\mu_0)$ and $\Ent_\infty(\mu_1)$ are finite.

In the case $N<\infty$ the argument goes as follows. For $\pi$-a.e. $\gamma\in\mathcal{A}$, we have that $d(\gamma_0,\gamma_1)>0$ and thus $\tau_{K,N}^{(t)}(d(\gamma_0,\gamma_1))>0$. Thus,
\begin{align}0<\int \tau_{K,N}^{(1-t)}(d(\gamma_0,\gamma_1))\rho_0(\gamma_0)+\tau_{K,N}^{(t)}(d(\gamma_0,\gamma_1))\rho_1(\gamma_1)\d\pi(\gamma)\le \Ent_N(\mu_t)=0
\end{align}
giving the contradiction.
\end{proof}
\begin{proposition}\label{LOCALISATION}
Let $(X,d,\m)$ be a metric measure space. Then $(X,d,\m)$ is very strict $CD(K,N)$ -space, if and only if for all absolutely continuous measures $\mu_0,\mu_1\in\P_2^{ac}(X)$ with bounded support, there exists an optimal plan $\pi\in\og(\mu_0,\mu_1)$, with $\mu_t\coloneqq(e_t)_\#\pi\in\P_2^{ac}(X)$,  for which the following two conditions hold:
\begin{enumerate}[(i)]
\item\label{cond1} For all $t\in(0,1)$, there exists a Borel map $T_t\colon X\to \geo(X)$ for which $\pi=(T_t)_\#\mu_t$, and $e_t\circ T_t=\id$.
\item\label{cond2} If $N=\infty$, then for every $t_1<t_2<t_3$, the inequality
\begin{align}\label{locconv}\log\rho_{t_2}(\gamma_{t_2})&\le \frac{(t_3-t_2)}{(t_3-t_1)}\log\rho_{t_1}(\gamma_{t_1})+\frac{(t_2-t_1)}{(t_3-t_1)}\log\rho_{t_3}(\gamma_{t_3})
\\&\phantom{=}-\frac{K}2\frac{(t_3-t_2)}{(t_3-t_1)}\frac{(t_2-t_1)}{(t_3-t_1)}d^2(\gamma_{t_1},\gamma_{t_3})
\end{align}
holds for $\pi$-almost every $\gamma$, where $\rho_t$ is the density of $\mu_t$ with respect to the reference measure $\m$.
\\\\
If $N<\infty$, then for every $t_1<t_2<t_3$, the inequality
\begin{align}\label{locconv2}\rho^{-\frac1N}(\gamma_{t_2})&\ge \tau_{K,N}^{\frac{(t_3-t_2)}{(t_3-t_1)}}(d(\gamma_{t_1},\gamma_{t_3}))\rho^{-\frac1N}_{t_1}(\gamma_{t_1})+\tau_{K,N}^{\frac{(t_2-t_1)}{(t_3-t_1)}}(d(\gamma_{t_1},\gamma_{t_3}))\rho^{-\frac1N}_{t_3}(\gamma_{t_3})
\end{align}
holds for $\pi$-almost every $\gamma$.
\end{enumerate}

Moreover, if $\pi$ is the plan given by the definition of very strict $CD(K,N)$ -space, then for $\pi$-almost every $\gamma$, the inequality \eqref{locconv} ($N=\infty$) or \eqref{locconv2} ($N<\infty$) holds for $\mathcal{L}^3$-almost every $(t_1,t_2,t_3)\in[0,1]$ with $t_1<t_2<t_3$.

\end{proposition}

\begin{remark} If we remove in Definition \ref{Sturm vsCD} the assumption of the boundedness of the supports of $\mu_0$ and $\mu_1$, we may remove it also from Proposition \ref{LOCALISATION}.
\end{remark}



\begin{proof}
We will prove only the case $N=\infty$. The proof of the finite dimensional case is the same with obvious modifications. Let $(X,d,\m)$ be a very strict $CD(K,\infty)$ -space. Let $\mu_0,\mu_1\in\P_2^{ac}(X)$, and let  $\pi$ be the optimal plan given by the definition of very strict $CD(K,\infty)$ -space. We will prove that the conditions \eqref{cond1} and \eqref{cond2} hold for $\pi$. By Lemma \ref{imabs} we have that $\mu_t$ is absolutely continuous with respect to $\m$. Moreover, the plan $(\restr_0^t)_\#\pi\in\og(\mu_0,\mu_t)$ is such as in the definition of very strict $CD(K,\infty)$. Thus, by Theorem \ref{exmap}, it is induced by a map $T$ from the intermediate measure $\mu_t$. Hence we have that $\pi=(S\circ e_0\circ T)_\#\mu_t\eqqcolon (T_t)_\#\mu_t$, where $S$ is the map given by Theorem \ref{exmap} for which $\pi=S_\#\mu_0$, proving the claim \eqref{cond1}.

For \eqref{cond2}, suppose to the contrary, that there exist $t_1,t_2,t_3\in[0,1]$, $t_1<t_2<t_3$, and a set $\mathcal{A}\subset \geo(X)$ with $\pi(\mathcal{A})>0$, so that the inequality \eqref{locconv} does not hold for any $\gamma\in \mathcal{A}$. Define $\tilde\pi\coloneqq \pi\lvert_\mathcal{A}$, and further define $\tilde\mu_t\coloneqq (e_t)_\#\tilde\pi=\tilde\rho_t\m$, for all $t\in[0,1]$.  Writing $\mathcal A$ as union
\[\mathcal A=\bigcup_{i\in\N} \{\gamma\in\mathcal A: \max_{j=1,2,3}{\rho_{t_j}(\gamma_{t_j})}\le i\},\] 
we may assume that $\tilde\rho_{t_j}$ is bounded from above, and so in particular that $(\tilde\rho_{t_j}\log \tilde\rho_{t_j})_+$ is integrable for $j\in\{1,2,3\}$. Let $\{\pi^t_x\}$ be the disintegration of $\pi$ with respect to the evaluation map $e_t$. Then we have, for all non-negative Borel functions $f\colon X\to \R$, that
\begin{align}\int_X f(x)\tilde\rho_t(x)\d\m(x)&=\int_{\geo(X)} f(\gamma_t)\chi_\mathcal{A}(\gamma)\d\pi(\gamma)
\\&=\int_X\int_{\geo(X)} f(\gamma_t)\chi_\mathcal{A}(\gamma)\d\pi^t_x(\gamma)\d\mu_t(x)
\\&=\int_X f(x)\int_{\geo(X)}\chi_\mathcal{A}(\gamma)\d\pi^t_x(\gamma)\d\mu_t(x)\\&=\int_X f(x)\left(\int_{\geo(X)}\chi_\mathcal{A}(\gamma)\d\pi^t_x(\gamma)\right)\rho_t(x)\d\m(x),
\end{align}
where $\rho_t$ is the density of $\mu_t\coloneqq (e_t)_\#\pi$ with respect to the reference measure $\m$. Thus $\tilde\rho_t(x)=\chi_\mathcal{A}(T_t(x))\rho_t(x)$ for $\m$-almost every $x\in X$. In particular, we have that \begin{align}\tilde\rho_t(\gamma_t)=\chi_\mathcal{A}(T_t(\gamma_t))\rho_t(\gamma_t)=\chi_\mathcal{A}(\gamma)\rho_t(\gamma_t), 
\end{align}
for $\pi$-almost every $\gamma$, and for all $t\in\{t_1,t_2,t_3\}$. Hence, we get
\begin{align}
\ent{\tilde\mu_{t_2}}&=\int_X\tilde\rho_{t_2}\log\tilde\rho_{t_2}\d\m=\int_X\log\tilde\rho_{t_2}\d\tilde\mu_{t_2}=\int_{\geo(X)}\log \tilde\rho_{t_2}(\gamma_{t_2})\d\tilde\pi
\\&=\int_\mathcal{A}\log\rho_{t_2}(\gamma_{t_2})\d\pi
\\&>\frac{(t_3-t_2)}{(t_3-t_1)}\int_\mathcal{A}\log\rho_{t_1}(\gamma_{t_1})\d\pi+\frac{(t_2-t_1)}{(t_3-t_1)}\int_\mathcal{A}\log\rho_{t_3}(\gamma_{t_3})\d\pi
\\&\phantom{\ge,}-\frac{K}2\frac{(t_3-t_2)}{(t_3-t_1)}\frac{(t_2-t_1)}{(t_3-t_1)}\int_\mathcal{A}d^2(\gamma_{t_1},\gamma_{t_3})\d\pi
\\&=\frac{(t_3-t_2)}{(t_3-t_1)}\ent{\tilde\mu_{t_1}}+\frac{(t_2-t_1)}{(t_3-t_1)}\ent{\tilde\mu_{t_3}}
\\&\phantom{=,}-\frac{K}2\frac{(t_3-t_2)}{(t_3-t_1)}\frac{(t_2-t_1)}{(t_3-t_1)}W_2^2(\tilde\mu_{t_1},\tilde\mu_{t_3}),
\end{align}
which contradicts the assumption of $\pi$ being the plan given by the definition of very strict $CD(K,\infty)$ -space. Hence \eqref{cond2} holds.

For the other direction, suppose that $\pi\in\og(\mu_0,\mu_1)$ is such that conditions \eqref{cond1} and \eqref{cond2} hold. Let $F\colon\geo(X)\to \R$ be a bounded non-negative Borel function for which $\int F\d\pi=1$, and let $t_1,t_2,t_3\in[0,1]$, $t_1<t_2<t_3$. Denote $\mu_t^F\coloneqq (e_t)_\#F\pi$. As previously, by \eqref{cond1}, we get that
\[\rho^F_t(x)\coloneqq F(T_t(x))\rho_t(x),\]
is the density of $\mu_t^F$ with respect to $\m$. Here $\rho_t$ is the density of $\mu_t$ with respect to the reference measure $\m$. In particular, we have that along geodesics the density is, up to a multiplicative constant, the same as the original density. More precisely, we have
\[\rho^F_t(\gamma_t)=F(\gamma)\rho_t(\gamma_t),\]
for $\pi$-almost every $\gamma$, and for every $t\in\{t_1,t_2,t_3\}$. Thus, by \eqref{cond2} we have that
\begin{align}\int\rho_{t_2}^F\log\rho_{t_2}^F\d\m&=\int\log\rho^F_{t_2}\d\mu_{t_2}^F=\int\log\rho_{t_2}^F(\gamma_{t_2})F(\gamma)\d\pi
\\&=\int\log\rho_{t_2}(\gamma_{t_2})F(\gamma)\d\pi+\int\log F(\gamma)F(\gamma)\d\pi
\\&\le\frac{(t_3-t_2)}{(t_3-t_1)}\int\log\rho_{t_1}(\gamma_{t_1})F(\gamma)\d\pi+\frac{(t_2-t_1)}{(t_3-t_1)}\int\log\rho_{t_3}(\gamma_{t_3})F(\gamma)\d\pi
\\&\phantom{=}-\frac{K}2\frac{(t_3-t_2)}{(t_3-t_1)}\frac{(t_2-t_1)}{(t_3-t_1)}\int d^2(\gamma_{t_1},\gamma_{t_3})F(\gamma)\d\pi+\int\log F(\gamma)F(\gamma)\d\pi
\\ &=\frac{(t_3-t_2)}{(t_3-t_1)}\int\rho_{t_1}^F\log\rho^F_{t_1}\d\mu_{t_1}^F+\frac{(t_2-t_1)}{(t_3-t_1)}\int\rho_{t_3}^F\log\rho_{t_3}^F\d\mu_{t_3}^F
\\&\phantom{=}-\frac{K}2\frac{(t_3-t_2)}{(t_3-t_1)}\frac{(t_2-t_1)}{(t_3-t_1)}W_2^2(\mu^F_{t_1},\mu_{t_3}^F),
\end{align}
giving the claim.

For the last claim, define for all $\gamma\in\geo(X)$ 
\[I_\gamma\coloneqq \{(t_1,t_2,t_3)\in J: \eqref{cond2}\mathrm{\ fails\ along\ } \gamma\mathrm{\ at\ }(t
_1,t_2,t_3)\},\]
where $J\coloneqq \{(t_1,t_2,t_3)\in [0,1]: t_1<t_2<t_3\}$, and the set 
\[\mathcal{I}\coloneqq \bigcup_\gamma \{\gamma\}\times I_\gamma.\]
Then by \eqref{cond2}
\begin{align}
0&=\int_J \pi(\{\gamma: t\in I_\gamma\})\d \L^3(t)=\int \chi_{\mathcal{I}}\d(\pi\otimes \L^3)
\\&=\int \L^3(I_\gamma)\d\pi(\gamma).
\end{align}
Hence $I_\gamma$ has Lebesgue measure zero for $\pi$-almost every $\gamma\in\geo(X)$.
\end{proof}

\begin{theorem}\label{equiviCD} Let $(X,d,\m)$ be a metric measure space.  Then the following are equivalent:
\begin{enumerate}[(i)]
\item \label{vsCD}The space $(X,d,\m)$ is a very strict $CD(K,N)$ -space (see Definition \ref{Sturm vsCD}). 
\item \label{VillanivsCD}The space $(X,d,\m)$ is a very strict $CD(K,N)$ -space in the spirit of $Lott--Villani$ (see Definition \ref{Lott--Villani vsCD}).
\end{enumerate}
\end{theorem}

\begin{proof}
Clearly condition \eqref{VillanivsCD} implies condition \eqref{vsCD}. For the other implication, assume that $\mu_0,\mu_1\in\P_2^{ac}(X)$, and $\pi\in\og(\mu_0,\mu_1)$ given by the definition of very strict $CD(K,N)$ -condition. Let $U\in\mathcal{DC}_N$, and $F\colon \geo(X)\to \R$ non-negative, bounded Borel function with $\int F\d\pi=1$. 

We first prove the claim, when $N<\infty$. Define $u(s)\coloneqq s^NU(s^{-N})$. Then $u$ is a decreasing and convex function, since $U\in\mathcal{DC}_N$. Hence, by Theorem \ref{LOCALISATION} condition \eqref{cond2},
\begin{align}
U(\mu_{t_2}^F)&=\int U\circ\rho^F_{t_2}\d\m= \int u((\rho^F_{t_2})^{-\frac1N})\rho^F_{t_2}\d\m =\int u((\rho^F_{t_2})^{-\frac1N}(\gamma_{t_2}))F(\gamma)\d\pi
\\ &= \int u(F(\gamma)\rho^{-\frac1N}_{t_2}(\gamma_{t_2}))F(\gamma)\d\pi 
\\ &\le \int u(F(\gamma)(\tau_{K,N}^\ymt(d(\gamma_{t_1},\gamma_{t_3}))\rho^{-\frac1N}_{t_1}(\gamma_{t_1})+\tau_{K,N}^\yt(d(\gamma_{t_1},\gamma_{t_3}))\rho_{t_3}^{-\frac1N}(\gamma_{t_3})))F(\gamma)\d\pi
\\&\le\ymt\int u(F(\gamma)\frac{(t_3-t_1)}{(t_3-t_2)}\tau_{K,N}^\ymt(d(\gamma_{t_1},\gamma_{t_3}))\rho^{-\frac1N}_{t_1}(\gamma_{t_1}))F(\gamma)\d\pi
\\&+\yt\int u(F(\gamma)\frac{(t_3-t_1)}{(t_2-t_1)}\tau_{K,N}^\yt(d(\gamma_{t_1},\gamma_{t_3}))\rho_{t_3}^{-\frac1N}(\gamma_{t_3}))F(\gamma)\d\pi
\\&= \frac{(t_3-t_2)}{(t_3-t_1)}U_{\pi,\m}^{\beta^{\frac{(t_3-t_2)}{(t_3-t_1)}}_{(K,N)}}(\mu^F_{t_1})+\frac{(t_2-t_1)}{(t_3-t_1)}U_{\pi^{-1},\m}^{\beta^{\frac{(t_2-t_1)}{(t_3-t_1)}}_{(K,N)}}(\mu^F_{t_3}),
\end{align}
giving the claim.

If $N=\infty$, we have that the function $u\colon s\mapsto e^sU(e^{-s})$ is convex and decreasing by assumption. Hence, by Proposition \ref{LOCALISATION} condition \eqref{cond2}
\begin{align} U(\mu_{t_2}^F)&= \int u(-\log(F(\gamma)\rho_{t_2}(\gamma_{t_2})))F(\gamma)\d\pi
\\ &\le\ymt\int u\left(-\log\left(\frac{F(\gamma)\rho_{t_1}(\gamma_{t_1})}{\beta_{\ymt}(\gamma_0,\gamma_1)}\right)\right)F(\gamma)\d\pi
\\ &\phantom{\le}+\yt\int u\left(-\log\left(\frac{F(\gamma)\rho_{t_3}(\gamma_{t_3})}{\beta_{\yt}(\gamma_0,\gamma_1)}\right)\right)F(\gamma)\d\pi
\\&=\ymt U_{\pi,\m}^{\beta_{(K,N)}^{\ymt}}(\mu^F_{t_1})+\yt U_{\pi^{-1},\m}^{\beta_{(K,N)}^{\yt}}(\mu^F_{t_3}),
\end{align} 
which completes the proof.
\end{proof}

Recall, that in our definition of very strict $CD(K,N)$-spaces, we only require the convexity of the entropy to hold for the critical exponent $N$, opposed to the definition of general $CD(K,N)$ -spaces. Therefore, the following immediate corollary is a non-trivial fact in this setting.
\begin{corollary}A metric measure space satisfying very strict $CD(K,N) (CD^*(K,N))$ -condition, satisfies very strict $CD(K,N') (CD^*(K,N'))$ -condition for any $N'>N$.\end{corollary}
\bibliographystyle{amsplain}
\bibliography{EquivsCDRef}
\end{document}